\documentclass[11pt,a4paper]{article}
\usepackage{amsfonts,latexsym,graphicx}

\newcommand{\G}{\Gamma}

\newcommand{\Z}{\mathbb{Z}}

\newcommand{\mod}{\mbox{mod}}
\newcommand{\rank}{\mbox{rank}}
\newcommand{\tr}{\mbox{tr}}

\title{The spectra of the spherical and euclidean triangle groups}
\author{Mark Harmer\\
Department of Mathematics\\
Australian National University\\
Canberra\\
email:mark.harmer@maths.anu.edu.au}
\date{}

\begin{document}

\maketitle 

\begin{abstract}
We derive the spectrum of the Laplace-Beltrami operator on the quotient orbifold of the non hyperbolic triangle groups. 
\end{abstract}

\section{Introduction}
We define the triangle groups
\begin{equation}\label{pres}
\G (p,q,r) = \left\langle \alpha , \beta\, |\, \alpha^p = \beta^q = (\alpha\beta)^r=e \right\rangle 
\end{equation}
where $p$, $q$ and $r$ are integers strictly greater than one. Those triangle groups for which
$$
\frac{1}{p} + \frac{1}{q} + \frac{1}{r} < 1 
$$
can be realised as discrete groups acting on the hyperbolic plane. The problem of finding the spectrum of the Laplace-Beltrami operator on the resulting quotient is, despite much effort, unsolved. The Selberg zeta function \cite{Sel}, the non trivial zeroes of which correspond to the eigenvalues of the Laplacian, may be used  to deduce some properties of the spectrum; however, it seems unlikely it may be used to find points on the spectrum. Currently, the only way to find eigenvalues is to numerically calculate the eigenmodes for the selected quotient space \cite{Hej1} (although Maass \cite{Maa}, who was first to consider this problem, was able to deduce some eigenvalues for a particular non co-compact group). \\
Here we will consider the non hyperbolic triangle groups, ie. those groups for which the above sum is greater than or equal to one. For the case of the above sum strictly greater than one the triangle group may be realised as a discrete group acting on the sphere, we refer to the triangle group as a spherical triangle group. If the sum is equal to one the triangle group may be realised as a discrete group acting on the plane and we refer to the group as an euclidean triangle group. It is perhaps not very surprising that for all of these non hyperbolic triangle groups the spectrum of the Laplacian on the resulting quotient space may be calculated explicitly. \\
In principle the Selberg zeta function associated to these triangle groups may be used to find the spectra. Our approach, we believe, emphasises the geometry of the problem. It appears that for each non hyperbolic triangle group we may easily calculate the spectrum of a normal torsion free subgroup (choosing the torsion free subgroup maximal is advisable). This subgroup is only non trivial for the euclidean triangle groups. It is then an easy step to find the spectrum of the original group by calculating the character of the representation of the quotient. 
\section{Enumeration of the non hyperbolic triangle groups}
We consider first the spherical triangle groups. Assuming $p\le q\le r$ we see that $p$ must be $2$ so that $q$ is $2$ or $3$. For $q=2$ we get the infinite family $\left\{ \G (2,2,n) \right\}_{n=2}$. For $q=3$  there are only three possibilities giving us the following spherical triangle groups
$$
\left\{ \G (2,2,n) \right\}^{<\infty}_{n=2} \, , \quad \G (2,3,3) \, , \quad \G(2,3,4) \, , \quad \G (2,3,5) \, .
$$
For the euclidean groups we again take $p\le q\le r$ so that $p$ must be $2$ or $3$. In the first case $q\le 4$ which gives $(q,r)\in\left\{ (2,\infty), (3,6), (4,4) \right\}$. In the second case $q\le 3$ which gives $(q,r)=(3,3)$. This gives the euclidean triangle groups
$$
\G (2,2,\infty) \, , \quad \G (2,3,6) \, , \quad \G(2,4,4) \, , \quad \G (3,3,3) \, .
$$
We will not consider the non co-compact case here apart from to observe that in the limit $n\to\infty$, $\G (2,2,n) \to \G (2,2,\infty)$ goes to the non co-compact euclidean triangle group. If we scale the radius of the sphere so that in the limit $n\to\infty$ we approach the plane it is an easy exercise to see that the spectra of $\G (2,2,n)$ accumulate to the continuous spectrum of $\G (2,2,\infty)$.
\section{The spherical triangle groups}
For the spherical triangle groups the torsion free subgroup is trivial, ie. the identity. The eigenfunctions invariant with respect to this group are of course the eigenfunctions on the sphere which are written in terms of associated Legendre functions
$$
\psi_{l,m} \left( \theta ,\varphi \right) = e^{im\varphi} P^m_l \left( \cos \theta \right)
$$
where $l\in\left\{ 0,1,\ldots\right\}$ and $m\in\left\{ -l, \ldots , l \right\}$. The corresponding eigenvalue is
$$
\lambda_l = l (l+1) \, .
$$
In this section we denote by $H_e (\lambda_l)$ the $2l+1$-dimensional eigenspace on the sphere corresponding to the eigenvalue $\lambda_l = l (l+1)$. \\
We will see that in the case of the spherical triangle groups the description of the spectrum reduces to studying the character of the representation of the particular triangle group on the space of eigenfunctions on the sphere.
\subsection{$\G (2,2,n)$}
The triangle groups $\G (2,2,n)$ are of course the dihedral groups $\Delta_n$. For this family we may explicitly calculate the eigenfunctions in terms of Legendre functions $\psi_{l,m}$ using the fact that the period $m$ should be a multiple of the order $n$. Instead we use the representation theory of $\G (2,2,n)$ as this is the approach we use for the other spherical groups and to a lesser degree for the euclidean groups. \\
We consider the projection 
\begin{equation}\label{projs}
P(\lambda_l)\,\psi= \frac{1}{\left| \G(2,2,n)\right|} \sum_{\gamma\in \G(2,2,n)} \psi \left( \gamma\, z \right)
\end{equation}
in $H_e (\lambda_l)$. Since we sum over $\G (2,2,n)$ this eigenfunction will clearly be invariant with respect to $\G (2,2,n)$, ie. $P(\lambda_l)$ is the identity on the eigenspace invariant with respect to $\G (2,2,n)$ which is certainly contained in $H_e (\lambda_l)$. This means that the multiplicity function $\mu_{\G} (\lambda_l)$ of the eigenvalue $\lambda_l$ for the group $\G (2,2,n)$ is 
$$
\mu_{\G} (\lambda_l) = \rank \left( P (\lambda_l) \right) = \tr \left( P (\lambda_l) \right) \, .
$$
For $\gamma\in\G(2,2,n)\subset SU(2)$ we can write $\gamma$ acting on $H_e (\lambda_l)$ as a matrix in terms of the basis $\{ \psi_{l,m}\}^{l}_{m=-l}$ of $H_e (\lambda_l)$. The components of this matrix can be found explicitly using the matrix elements of the regular representation of $SU (2)$ \cite{Vil}. However, all we need is the trace, or character of the representation, which (\cite{Vil}, pg 170) is constant under conjugacy so depends only on the angle $\chi$ of rotation due to $\gamma$ through
\begin{equation}\label{char}
\tr \left( \gamma (\lambda_l) \right) = \frac{\sin \left(l+\frac{1}{2}\right)\chi}{\sin \frac{\chi}{2}} \, .
\end{equation}
We use $\gamma (\lambda_l)$ to indicate that we are considering the action of the rotation on functions in the space $H_e (\lambda_l)$. Consequently 
\begin{eqnarray*}
\mu_{\G} (\lambda_l) & = & \tr \left( P (\lambda_l) )\right) = \frac{1}{\left| \G(2,2,n)\right|} \sum_{\gamma_i\in \G(2,2,n)} \frac{\sin \left(l+\frac{1}{2}\right)\chi_i}{\sin \frac{\chi_i}{2}} \\
& = & \frac{1}{2n} \left( 2l+1 + (-1)^l \, n + \sum^{n-1}_{s=1}
\frac{\sin \left(2l+1\right)\frac{s\pi}{n}}{\sin \frac{s\pi}{n}} \right) \, .
\end{eqnarray*}
Here we use the fact that $\G (2,2,n)$ has $2n$ elements; an element of order $n$ plus $n$ rotations of order two. \\
We can simplify this expression for the multiplicity by using the formula of Eisenstein \cite{Eis}
\begin{equation}\label{ETS}
\left(\left( \frac{l}{n} \right)\right) = -\frac{1}{2n} \sum^{n-1}_{s=1}
\sin \left( \frac{2l s\pi}{n} \right) \cot \left( \frac{s\pi}{n} \right)
\end{equation}
where
$$
\left(\left( x \right)\right) = \left\{ \begin{array}{cl}
x - \lfloor x \rfloor - \frac{1}{2} & x\not\in\Z \\
0 & x\in\Z 
\end{array} \right.
$$
and $\lfloor x \rfloor$ is the largest integer less than $x$. Writing the trigonometric sum 
which appears in the multiplicity as 
$$
\sum^{n-1}_{s=1}
\frac{\sin \left(2l+1\right)\frac{s\pi}{n}}{\sin \frac{s\pi}{n}} = \sum^{n-1}_{s=1}
\left[  \sin \left( \frac{2l s\pi}{n} \right) \cot \left( \frac{s\pi}{n} \right) + \cos \left( \frac{2l s\pi}{n} \right) \right]
$$
and using (\ref{ETS}) we arrive at the following simpler expression for the multiplicity
$$
\mu_{\G} \left( \lambda_l \right) = \left\lfloor \frac{l}{n} \right\rfloor + \frac{1 + (-1)^{l}}{2} \, .
$$
We note that the formula for the multiplicity can also be simplified using the above mentioned description of the spectrum of the dihedral family in terms of associated Legendre functions. 
\subsection{$\G (2,3,3)$}
It is clear that $\G (2,3,3) = A_4$, the pure rotation symmetry group of the tetrahedron. We consider the projection  $P(\lambda_l)$ in $H_e (\lambda_l)$ defined as in (\ref{projs}) by a sum over $\G(2,3,3)$.
Again the multiplicity function will be the trace of this projection. Using the fact that $\G (2,3,3)$ has three elements of order two and four elements of order three, for a total of twelve elements, along with (\ref{char}) we evaluate the trace as
\begin{eqnarray*}
\mu_{\G} (\lambda_l) & = & \frac{1}{12} \left( 2l+1 + (-1)^l \, 3 + 4 \sum^{2}_{s=1}
\frac{\sin \left(2l+1\right)\frac{s\pi}{3}}{\sin \frac{s\pi}{3}} \right) \\
& = & 2 \left\lfloor \frac{l}{3} \right\rfloor + \frac{3 + (-1)^{l}-2l}{4} \, .
\end{eqnarray*}
In the last line we use (\ref{ETS}) to simplify. 
\subsection{$\G (2,3,4)$}
It is clear that $\G (2,3,4) = S_4$, the pure rotation symmetry group of the cube and octahedron. Analogous to (\ref{projs}) we define the projection $P(\lambda_l)$ onto the eigenspace invariant with respect to $\G (2,3,4)$ and calculate the trace using (\ref{char}). Using the fact that $\G (2,3,4)$ has six elements of order two, four elements of order three and three elements of order four for a total of twenty four elements we get the multiplicity function
\begin{eqnarray*}
\mu_{\G} (\lambda_l) & = & \frac{1}{24} \left( 2l+1 + (-1)^l \, 6 + 
4 \sum^{2}_{s=1} \frac{\sin \left(2l+1\right)\frac{s\pi}{3}}{\sin \frac{s\pi}{3}} + 
3 \sum^{3}_{s=1} \frac{\sin \left(2l+1\right)\frac{s\pi}{4}}{\sin \frac{s\pi}{4}} \right) \\
& = & \left\lfloor \frac{l}{3} \right\rfloor + \left\lfloor \frac{l}{4} \right\rfloor + \frac{3 + (-1)^{l}-2l}{4} \, .
\end{eqnarray*}
\subsection{$\G (2,3,5)$}
It is clear that $\G (2,3,5) = A_5$, the pure rotation symmetry group of the dodecahedron and icosahedron. Analogous to (\ref{projs}) we define the projection $P(\lambda_l)$ onto the eigenspace invariant with respect to $\G (2,3,5)$ and calculate the trace using (\ref{char}). Using the fact that $\G (2,3,5)$ has fifteen elements of order two, ten elements of order three and six elements of order five for a total of sixty elements we get the multiplicity function
\begin{eqnarray*}
\mu_{\G} (\lambda_l) & = & \frac{1}{60} \left( 2l+1 + (-1)^l \, 15 + 10 \sum^{2}_{s=1}
\frac{\sin \left(2l+1\right)\frac{s\pi}{3}}{\sin \frac{s\pi}{3}} + 6 \sum^{4}_{s=1}
\frac{\sin \left(2l+1\right)\frac{s\pi}{5}}{\sin \frac{s\pi}{5}} \right) \\
& = & \left\lfloor \frac{l}{3} \right\rfloor + \left\lfloor \frac{l}{5} \right\rfloor + \frac{3 + (-1)^{l}-2l}{4} \, . 
\end{eqnarray*}
\section{The euclidean triangle groups}
Here the normal torsion free subgroups are a little more interesting. We will see that the spectrum can be most simply described by choosing a maximal subgroup. In this section we ignore the $\lambda=0$ eigenvalue which will always have multiplicity one.
\subsection{$\G (2,3,6)$}
We know that this group generates an action on the plane with fundamental domain which is an equilateral triangle---we assume that the triangle has side length $4\pi /3$. Using the presentation (\ref{pres}) we consider the elements
$$
\tau = (\alpha\beta)^3 \alpha \, , \quad \sigma = (\alpha\beta)^2 \beta \, .
$$
In terms of the group action on the plane these elements are both translations so that the group generated by them,
$$
T(2,3,6) = \left\langle \tau , \tau^{-1} , \sigma , \sigma^{-1} \right\rangle \, ,
$$
is torsion free. Furthermore,
\begin{eqnarray*}
\alpha \tau \alpha^{-1} = \tau^{-1} \, , & \quad & \beta \tau \beta^{-1} = \tau^{-1} \sigma \\
\alpha \sigma \alpha^{-1} = \sigma^{-1} \, , & \quad & \beta \sigma \beta^{-1} = \tau^{-1} 
\end{eqnarray*}
so that $T(2,3,6)\lhd\G (2,3,6)$. The quotient may be identified with
$$
\G (2,3,6)/T(2,3,6) = S(2,3,6) = \left\{ e, \gamma , \gamma^2 ,\ldots , \gamma^5 \right\} \, ,
$$
where $\gamma=\alpha\beta$. This follows since $T(2,3,6)\cap S(2,3,6) = \{e\}$ and the join
$$
T(2,3,6) \vee S(2,3,6) = \G (2,3,6) \, ,
$$
which in turn follows if we observe that a fundamental domain of $T(2,3,6)$ is the hexagon while $S(2,3,6)$ permutes the fundamental domains of $\G (2,3,6)$ in the hexagon. \\
We first solve the spectral problem for $T(2,3,6)$: we find eigenfunctions of the Laplacian
$$
-\left( \frac{\partial^2 \psi}{\partial x^2} + \frac{\partial^2 \psi}{\partial y^2} \right) = \lambda \psi
$$
on the hexagonal lattice, ie. subject to 
\begin{eqnarray*}
&& \psi (\tau\, z) = \psi \left( x+\frac{4\pi}{\sqrt{3}}, y \right) = \psi (x,y) \\
&& \psi (\sigma\, z) = \psi \left( x+\frac{2\pi}{\sqrt{3}}, y+ 2\pi \right) = \psi (x,y) \, .
\end{eqnarray*}
By a standard result we obtain a complete list of the eigenfunctions in the form $\exp \left( i (a x + b y) \right)$ where $(a,b)$ lies in the lattice dual to $T(2,3,6)$. Explicitly, the eigenfunctions are
$$
\psi_{m,n} (x,y) = \exp \frac{i}{2} \left( \sqrt{3} m x + (2n+m) y \right)
$$
with corresponding eigenvalue
$$
\lambda_{m,n} = m^2 + n^2 + mn 
$$
with $(m,n)\in\Z^2$. Using some number theory \cite{Hir} we immediately have the multiplicity function for $T(2,3,6)$
\begin{eqnarray*}
\mu_T (\lambda) & = & \# \left\{ (m,n)\in\Z^2 \, \left|\, \lambda = m^2 + n^2 + mn \right. \right\} \\
& = & 6 \left( d_{1,3} \left( \lambda \right) - d_{2,3} \left( \lambda \right) \right) \, .
\end{eqnarray*}
Here $d_{r,s} (N)$ is the number of divisors of $N$ which are equivalent to $r$ mod $s$. \\
Let us denote by $H_T (\lambda)$ the $\mu_T (\lambda )$-dimensional eigenspace of $T (2,3,6)$ corresponding to eigenvalue $\lambda$. We consider the projection 
\begin{equation}\label{projp}
P(\lambda)\, \psi = \frac{1}{\left| S(2,3,6)\right|} \sum_{\gamma\in S(2,3,6)} \psi \left( \gamma\, z \right) 
\end{equation}
in $H_T (\lambda)$ where the sum is now over the quotient group. Using $\G (2,3,6) =T(2,3,6)\vee S(2,3,6)$ and $T(2,3,6)\lhd\G(2,3,6)$ it is not difficult to see that $P(\lambda)\,\psi$ is an eigenfunction which is invariant with respect to $\G (2,3,6)$. As above $P (\lambda)$ is unity on the eigenspace of $\G (2,3,6)$ which is certainly a subspace of $H_T (\lambda)$ so that
$\tr \left( P(\lambda ) \right) = \rank \left( P (\lambda) \right) = \mu_{\G} (\lambda)$
is the multiplicity function for the group $\G (2,3,6)$. \\
Calculating the trace of $P(\lambda)$ is a little easier than above: we consider the action of the $2\pi /6$ rotation
$$
\gamma \, \left( \begin{array}{c}
x \\ y
\end{array} \right) = \left( \begin{array}{cc}
\frac{1}{2} & -\frac{\sqrt{3}}{2} \\
\frac{\sqrt{3}}{2} & \frac{1}{2}
\end{array} \right) \left( \begin{array}{c}
x \\ y
\end{array} \right)
$$
on the basis $\psi_{m,n}$ of $H_T (\lambda)$. It is easy to see that $\gamma^i$, for $i\not\equiv 0\, \mod\, (6)$, permutes the basis elements $\psi_{m,n}$ so as a matrix it will only have zeroes on the diagonal. Consequently, the trace of $P (\lambda)$ will just be the trace of the identity divided by the order of $S(2,3,6)$ or
$$
\mu _{\G} (\lambda) = \frac{1}{6} \mu_T (\lambda) = d_{1,3} \left( \lambda \right) - d_{2,3} \left( \lambda \right) \, .
$$
\subsection{$\G (2,4,4)$}
Here the fundamental domain is a square which we assume to have side length $\pi$. Using the presentation (\ref{pres}) we consider the torsion free subgroup $T (2,4,4)$ generated by the translations
$$
\tau = \beta^2 \alpha \, , \quad \sigma = (\alpha\beta)^2 \alpha \, .
$$
Since
\begin{eqnarray*}
\alpha \tau \alpha^{-1} = \tau^{-1} \, , & \quad & \beta \tau \beta^{-1} = \sigma \\
\alpha \sigma \alpha^{-1} = \sigma^{-1} \, , & \quad & \beta \sigma \beta^{-1} = \tau^{-1} 
\end{eqnarray*}
$T (2,4,4)$ is a normal subgroup with, for the same reasons as above, quotient
$$
\G (2,4,4)/T(2,4,4) = S(2,4,4) = \left\{ e, \gamma , \gamma^2 , \gamma^3 \right\} \, ,
$$
where again $\gamma=\alpha\beta$. \\
The spectral problem for $T (2,4,4)$ amounts to finding eigenfunctions on the square lattice, ie. subject to
\begin{eqnarray*}
&& \psi (\tau\, z ) = \psi \left( x+2\pi, y \right) = \psi (x,y) \\
&& \psi (\sigma\, z) = \psi \left( x, y+2\pi \right) = \psi (x,y) \, .
\end{eqnarray*}
The eigenfunctions and eigenvalues for $T (2,4,4)$ are then
\begin{eqnarray*}
\psi_{m,n} (x,y) & = & \exp i \left( m x + n y \right) \\
\lambda_{m,n} & = & m^2 + n^2 
\end{eqnarray*}
with $(m,n)\in\Z^2$. Again \cite{CoLa} we immediately have the multiplicity function for $T(2,4,4)$
\begin{eqnarray*}
\mu_T (\lambda) & = & \# \left\{ (m,n)\in\Z^2 \, \left|\, \lambda = m^2 + n^2 \right. \right\} \\
& = & 4 \left( d_{1,4} \left( \lambda \right) - d_{3,4} \left( \lambda \right) \right) \, .
\end{eqnarray*}
Returning to the original triangle group we see that, as described above, the multiplicity function for $\G (2,4,4)$ can be found from the trace of the projection $P(\lambda)$ defined as in (\ref{projp}) by taking the sum over the quotient $S (2,4,4)$. Once again, the matrix of $\gamma^i$, for $i\not\equiv 0\, \mod\, (4)$, in the basis $\psi_{m,n}$ will only have zeroes on the diagonal. Consequently, the trace of $P (\lambda)$ will just be the trace of the identity divided by the order:
$$
\mu _{\G} (\lambda) = \frac{1}{4} \mu_T (\lambda) = d_{1,4} \left( \lambda \right) - d_{3,4} \left( \lambda \right) \, .
$$
\subsection{$\G (3,3,3)$}
The fundamental domain of $\G (3,3,3)$ consists of a rhombus with acute angle $\pi /3$. We assume the side of the rhombus has length $4\pi /3$. Using (\ref{pres}) we define the translations
$$
\tau = \beta^2 \alpha \, , \quad \sigma = \alpha\beta \alpha \, .
$$
The group $T (3,3,3)$ generated by these translations is normal in $\G (3,3,3)$
\begin{eqnarray*}
\alpha \tau \alpha^{-1} = \tau^{-1}\sigma \, , & \quad & \beta \tau \beta^{-1} = \tau^{-1} \sigma \\
\alpha \sigma \alpha^{-1} = \tau^{-1} \, , & \quad & \beta \sigma \beta^{-1} = \tau^{-1} 
\end{eqnarray*}
with quotient
$$
\G (3,3,3)/T(3,3,3) = S(3,3,3) = \left\{ e, \gamma , \gamma^2  \right\} \, .
$$
It is useful for us to observe that $T (3,3,3)$ and $T(2,3,6)$ are the same group ($\G (3,3,3)$ is a subgroup of $\G (2,3,6)$ and both have the same maximal subgroup generated by translations). Consequently, $T (3,3,3)$ has the same eigenfunctions, eigenvalues and in particular multiplicity function 
$$
\mu_T (\lambda) = 6 \left( d_{1,3} \left( \lambda \right) - d_{2,3} \left( \lambda \right) \right) \, .
$$
as $T (2,3,6)$. \\
The projection $P(\lambda)$ is a sum over $S (3,3,3)$ which is a subgroup of $S (2,3,6)$ and as we ascertained above the matrix elements for the terms of this sum $\gamma^i$, for $i\not\equiv 0\, \mod\, (3)$, will only have zeroes on the diagonal. Consequently, the multiplicity function for $\G (3,3,3)$ is 
$$
\mu _{\G} (\lambda) = \frac{1}{3} \mu_T (\lambda) = 2 \left( d_{1,3} \left( \lambda \right) - d_{2,3} \left( \lambda \right) \right) \, .
$$
It is not difficult to see that the spectral problem for the $\G (3,3,3)$ triangle group is equivalent to solving the two eigenvalue problems for the Laplacian on the equilateral triangle with Neumann and Dirichlet boundary conditions. This problem has been considered by many authors \cite{MCar, Pin1, Pin2, Pra} right back to Lam\'{e} \cite{Lam}. \\
It seems that the approach to this problem most similar to our approach is to be found in the papers by Pinsky. In \cite{Pin1, Pin2} the author, in effect, considers a subgroup of translations which is not maximal---it has fundamental domain consisting of nine copies of the fundamental domain of $\G (3,3,3)$. Consequently his description of the spectrum is a little more complicated (hence the comment above that it is beneficial to consider the maximal torsion free subgroup). The reason for this extra complication is that an important issue in all of the cited papers is to distinguish between Neumann and Dirichlet eignevalues.
\section{Summary}
In table \ref{multitab} we summarise the multiplicity functions for each of the (co-compact) non hyperbolic triangle groups.  Subject to the normalisations given above the spectrum is always a subset of the integers, in the case of the spherical triangle groups of the form $\lambda_l = l (l+1)$, $l\in\{0,1,\ldots\}$.
We also note that, in the case of the euclidean groups, $\mu (0) = 1$, the multiplicity at zero is always one. \\
\begin{table}[h]
\begin{center}\hspace*{-8mm}
\begin{tabular}{|c|c|}  \hline
$\G (2,2,n)$ & $\rule[-3mm]{0mm}{8mm} \mu \left( \lambda_l \right) = \left\lfloor \frac{l}{n} \right\rfloor + \frac{1 + (-1)^{l}}{2}$ \\ \hline
$\G (2,3,3)$ & $\rule[-3mm]{0mm}{8mm} \mu(\lambda_l) = 2 \left\lfloor \frac{l}{3} \right\rfloor + \frac{3 + (-1)^{l}-2l}{4}$ \\ \hline
$\G (2,3,4)$ & $\rule[-3mm]{0mm}{8mm} \mu(\lambda_l) = \left\lfloor \frac{l}{3} \right\rfloor + \left\lfloor \frac{l}{4} \right\rfloor + \frac{3 + (-1)^{l}-2l}{4}$ \\ \hline
$\G (2,3,5)$ & $\rule[-3mm]{0mm}{8mm} \mu(\lambda_l) = \left\lfloor \frac{l}{3} \right\rfloor + \left\lfloor \frac{l}{5} \right\rfloor + \frac{3 + (-1)^{l}-2l}{4}$ \\ \hline
$\G (2,3,6)$ & $\rule[-2mm]{0mm}{6mm} \mu(\lambda) = d_{1,3} \left( \lambda \right) - d_{2,3} \left( \lambda \right)$ \\ \hline
$\G (2,4,4)$ & $\rule[-2mm]{0mm}{6mm} \mu(\lambda) = d_{1,4} \left( \lambda \right) - d_{3,4} \left( \lambda \right)$ \\ \hline
$\G (3,3,3)$ & $\rule[-2mm]{0mm}{6mm} \mu(\lambda) = 2 \left( d_{1,3} \left( \lambda \right) - d_{2,3} \left( \lambda \right) \right)$ \\ \hline
\end{tabular}
\end{center}
\caption{Multiplicity functions of the spherical and euclidean triangle groups.}\label{multitab}
\end{table}
We do not attempt to explicitly calculate the eigenfunctions here apart from to state that, in the case of the euclidean groups, they can always be calculated from the eigenfunctions of the appropriate torus by applying the finite sum (\ref{projp}) over the associated quotient group. In the case of the spherical triangle groups we have noted that the eigenfunctions of the dihedral groups may be explicitly calculated in terms of associated Legendre functions. Then we may use the fact that $\G (2,2,2) \lhd \G (2,3,3) \lhd \G (2,3,4)$ and $\G(2,3,3) \subset \G (2,3,5)$ to write the eigenfunctions of $\G (2,3,3)$, $\G (2,3,4)$ and $\G (2,3,5)$ as finite sums of eigenfunctions of the preceding group in the list. \\
Of some independent interest are the high energy asymptotics of the counting function
$$
N_{\G} \left( \Lambda \right) = \sum^{\Lambda}_{\lambda} \mu_{\G} (\lambda) \, ,
$$
the so called Weyl asymptotics. In the case of the spherical triangle groups these can be calculated explicitly. To do this we start with the dihedral family $\G (2,2,n)$ for which we can explicitly calculate the counting function
\begin{eqnarray*}
N_{\G(2,2,n)} \left( L \right) & = & 1 + \left\lfloor \frac{L}{2} \right\rfloor + ( L-n+1) + \cdots + (L-kn+1) + \cdots \\
& = & \frac{(L+1)^2}{2n} + O (1) \, .
\end{eqnarray*}
Here we should recall that the eigenvalue is related to $L$ by $\Lambda = L (L+1)$. This implies that the oscillating term 
$$
\sum^{n-1}_{s=1} \frac{\sin \left(2l+1\right)\frac{s\pi}{n}}{\sin \frac{s\pi}{n}} \, ,
$$
which appears in the multiplicity function of each of the other spherical triangle groups, remains bounded when we sum over $l$. Consequently, the counting functions for the remaining spherical triangle groups have the following asymptotic form
\begin{eqnarray*}
N_{\G(2,3,3)} \left( L \right) & = & \frac{(L+1)^2}{12} + O (1) \\
N_{\G(2,3,4)} \left( L \right) & = & \frac{(L+1)^2}{24} + O (1) \\
N_{\G(2,3,5)} \left( L \right) & = & \frac{(L+1)^2}{60} + O (1) \, .
\end{eqnarray*}
The euclidean triangle groups are unfortunately not so simple. From geometric considerations it is clear that 
the first term in the asymptotic expansion will be $\frac{\pi}{4}\lambda$ for $\G (2,4,4)$ (the square lattice) and $\frac{2\pi}{\sqrt{3}|S|}\lambda$ for $\G (2,3,6)$ and $\G (3,3,3)$ (the hexagonal lattice) where $S$ is the quotient of the triangle group by its torsion free subgroup. The next term is the subject of a difficult open conjecture in analytic number theory \cite{Gro,Hardy}, suffice to say that we do not even know the order of this term (it is conjectured to be $O(\lambda^{\frac{1}{4}+\epsilon})$).
\section*{Acknowledgements}
While doing this work the author was supported by a New Zealand FRST postdoctoral fellowship. The author would like to acknowledge financial assistance from the Marsden fund as well as useful conversations with Prof. Gaven Martin, Dr Shaun Cooper and Dr Heng Huat Chan. The author would like to thank one of the referees for comments which led to a simplification of the multiplicity functions.

\end{document}